\documentclass[12pt]{article}

\usepackage{amsmath, amsthm, amscd, amsfonts, amssymb}

\newtheorem{theorem}{Theorem}[section]
\newtheorem{proposition}[theorem]{Proposition}
\newtheorem{corollary}[theorem]{Corollary}
\newtheorem{lemma}[theorem]{Lemma}
\theoremstyle{definition}
\newtheorem{definition}[theorem]{Definition}

\newtheorem{remark}[theorem]{Remark}

\def\cald{{\mathcal{D}}}
\def\calx{{\mathcal{X}}}

\def\calk{{\mathcal{K}}}

\def\calq{{\mathcal{Q}}}
\def\calx{{\mathcal{X}}}

\def\pf{\n{\bf Proof.} }
\def\vsp{\vspace*{1,5mm}\\ }
\def\bk{\bigskip }

\def\n{\noindent }
\def\dd{\displaystyle}

\def\barr{\begin{array}}
\def\earr{\end{array}}

\def\rr{{\mathbb{R}}}

\def\1{^{-1}}
\def\9{{\infty}}
\def\lbb{{\lambda}}

\def\wt{\widetilde}
\def\ov{\overline}
\def\vf{{\varphi}}

\def\ooo{{\Omega}}
\def\pp{{\partial}}
\def\vp{{\varepsilon}}

\def\ff{\forall }
\def\({\left(}
\def\){\right)}
\def\<{\left<}
\def\>{\right>}

\title{Boundary controllability of phase-transition region of~a~two-phase Stefan problem}
\author{Viorel Barbu\thanks{Octav Mayer Institute of Mathematics of the Romanian Academy and Al.I. Cuza University, Ia\c si, Romania.  Email: vb41@uaic.ro}}
\date{}

\begin{document}
\maketitle
\begin{abstract}
\n One proves that the moving interface  of a two-phase Stefan problem on $\ooo\subset\rr^d$, $d=1,2,3,$ is   controllable at the end time $T$ by a Neumann boundary controller $u$. The phase-transition region is a mushy region $\{\sigma^u_t;\ 0\le t\le T\}$ of a modified Stefan problem and the main result amounts to saying that, for each Lebesque measurable set $\ooo^*$ with positive measure, there is $u\in L^2((0,T)\times\pp\ooo)$ such that  $\ooo^*\subset\sigma^u_T.$ To this aim, one uses an optimal control approach combined with Carleman's inequality and the Kakutani fixed point theorem.
\end{abstract}

\section{Introduction}\label{s1}

The heat conduction in a bounded, open domain $\ooo\subset\rr^d,$ $d=1,2,3,$ with two-phase transition (ice-water) is described by the classical  Stefan problem
\begin{equation}\label{e1.1}
\barr{ll}
\dd\frac{\pp\theta}{\pp t} -k_1\Delta\theta=0\mbox{\ in }\calq_-=\{(t,x)\in\calq;\ \theta(t,x)<0\},\vsp
\dd\frac{\pp\theta}{\pp t} -k_2\Delta\theta=0\mbox{\ in }\calq_+=\{(t,x)\in\calq;\ \theta(t,x)>0\},\vsp
(k_2\nabla_x\theta^+(t,x)-k_1\nabla_x\theta^-(t,x))\cdot n_x(t,x)=-\rho\vsp\hfill\mbox{ on }S=\{(t,x)\in\calq,\ \theta^-(t,x)=0\},\\
\theta(0,x)=\theta_0(x),\ \ x\in\ooo.
\earr\end{equation}
Here $\theta=\theta(t,x)$ is the temperature in $(t,x)\in\calq=(0,T)\times\ooo$, $(n_t,n_x)$ is the outward normal to $S$ with $\theta^+=\theta$ in $\calq_+,$ $\theta^-=\theta$ in $\calq_-$. Moreover, $k_1,k_2$ are heat conductivities in solid and liquid regions and $\rho>0$ is the latent heat. If $\Gamma=\pp\ooo$ is the boundary of $\ooo$, which is assumed smooth (of class $C^1$, for instance), and a heat flux $u$ is applied on $\Gamma,$ then the Neumann boundary condition
\begin{equation}\label{e1.2}
\frac{\pp\theta}{\pp \nu}=u\mbox{\ \ on }\Sigma=(0,T)\times\Gamma\end{equation}should be added to \eqref{e1.1}. (Here $\frac{\pp\theta}{\pp\nu}=\nabla\theta\cdot\nu$ and $\nu$ is the outward normal to~$\Gamma$.)
In the following, $u\in L^1(\Sigma)$ is a {\it boundary controller}.

If we denote by $\beta:\rr\to\rr$ the enthalpy of system, that is,
\begin{equation}\label{e1.3}
\beta(r)=\left\{\barr{ll}
k_1r&\ff\,r<0,\vsp
0&\ff\,r\in[0,\rho),\vsp
k_2(r-\rho)&\ff\,r>\rho,\earr\right.
\end{equation}
one can rewrite \eqref{e1.1} as
\begin{equation}\label{e1.4}
\barr{ll}
\dd\frac{\pp}{\pp t}\ \beta^{-1}(\theta)-\Delta\theta=0&\mbox{ in }\cald'(\calq),\vsp
\theta(0)=\theta_0&\mbox{ in }\ooo,\vsp
\dd\frac{\pp\theta}{\pp\nu}=u&\mbox{ on }\Sigma.\earr
\end{equation}
Equivalently (see \cite{2}, \cite{3}, \cite{5}, \cite{11}),
\begin{equation}\label{e1.5}
\barr{ll}
\dd\frac{\pp y}{\pp t}-\Delta\beta(y)=0&\mbox{ in }\calq,\vsp
y(0)=y_0=\beta^{-1}(\theta_0)&\mbox{ in }\ooo,\vsp
\dd\frac{\pp}{\pp\nu}\,(\beta(y))=u&\mbox{ on }\Sigma,\earr
\end{equation}
where $\beta(y)=\theta.$

It turns out that (see, e.g., \cite{3}, \cite{11}) that, for each $u\in L^2(\Sigma)$ and $y_0\in L^2(\ooo),$ equation \eqref{e1.5} has a unique weak (variational) solution $y=u$. If the interface $S=S^u=\{(t,x);\ \theta(t,x)=0\}=\{(t,x);\ y^u(t,x)\in[0,8]\}$ is smooth, we may represent it as
$$S=\bigcup_{t\in[0,T]}\Gamma^u_t,\ \Gamma^u_t=\{x\in\ooo;\theta^u(t,x)=0\}=\{x\in\ooo;y^u(t,x)\in[0,\rho]\},$$
where $\Gamma^u_t$ is a smooth {\it moving boundary}.\newpage

The standard controllability problem of the moving boundary $\Gamma^u_t$ is the following:
\begin{quote}
	{\it Given a surface $\Gamma^*\subset\ooo$, find $u^*\in L^2(\Sigma)$ such that $\Gamma^{u^*}_T=\Gamma^*.$}
\end{quote}
This problem is very important in technology (continuous casting of steel, for instance) and was intensively studied in the last decades (see, e.g., \cite{12}, \cite{14}).  Most of the results refer to the numerical treatment of such controllability problem with a few exceptions including \cite{3}, \cite{3a}, \cite{6}, \cite{8}. We note, for instance, that such a real time feedback controller in $1-D$ is designed in \cite{4} is designed.  However, in general, the problem remained open.

On the other hand, it should be mentioned that the above classical Stefan model of phase transition is quite inappropriate to describe the real physical process which exhibits superheating and leads so to unstability of the phase-change surface $S$. (See, e.g., \cite{1}, \cite{9}.) In fact, a more realistic model is the melting process with three regions: the {\it solid} one,  $Q^t_-{=}\{x\in\ooo;\theta(t,x)<-\mu\}$, the {\it liquid}, $Q^t_+=\{x\in\ooo;\theta(t,x)>\mu\}$ and the {\it mushy} region $\sigma_t=\{x\in\ooo;$ $|\theta(t,x)|\le\mu\}$, where the temperature $\theta$ equals to melting temperature, while the substance (the ice or water in our case) is neither pure solid nor pure liquid. Thus, in reality, the interface $\Gamma^u_t$ is no more a smooth surface but the mushy region$\sigma_t=\sigma^u_t$ equivalently represented as
\begin{equation}
\label{e1.6}\sigma^u_t=\{x\in\ooo;-\mu\le y^u(t,x)\le\rho+\mu\}.
\end{equation}
Then the controllability problem for the moving boundary $\Gamma^u_t$ reduces to the more realistic situation of the controllability of mushy region $\sigma^u_t$. Namely, {\it given $T>0$ and $\ooo^*\subset\ooo$, find $u^*\in L^2(\Sigma)$ such that ${\rm dist}_H(\sigma^{u^*}_T,\ooo^*)=0.$} (${\rm dist}_H$ is the Hausdorff-Pompeiu distance.)

Here, we shall prove (Theorem \ref{t3.1}) a slightly weaker version of this problem. Namely,
\begin{quote}\it Given $\ooo^*\subset\ooo$, there is $u^*\in L^2(\Sigma)$ such that $\ooo^*\subset\sigma^{u^*}_T.$
	\end{quote}
As a matter of fact, this result will be proven for a regular modified version of \eqref{e1.1} which reduces to \eqref{e1.1} on pure liquid and solid regions. This model will be presented and analyzed in Section \ref{s2} and the regularity of its solution (H\"older continuity) necessary to give a meaning to this phase-transition model represents much of the substance of this work.\bk

\n{\bf Notations.} $\ooo$ is an open and bounded domain with smooth boundary $\Gamma$, $\calq=(0,T)\times\ooo,$ $\Sigma=(0,T)\times\Gamma;$ $L^p(\ooo)$, $L^p(\calq)$ and $L^p(\Sigma)$, $1\le p\le\9$, are the $L^p$-Lebesgue integrable functions on $\ooo$, $\calq$ and $\Sigma$, respectively, with the norms denoted $|\cdot|_p,$ $\|\cdot\|_{L^p(\calq)}$ and $\|\cdot\|_{L^p(\Sigma)}$. Denote by $W^1_p(\ooo)$, \mbox{$1\le p\le\9$,} the Sobolev space
$\left\{y\in L^p(\ooo);\frac{\pp y}{\pp x_1}\in L^p(\ooo),\ i=1,2,,...,d\right\}$ and by $V$ the   space $H^1(\ooo)=W^1_2(\ooo)=\left\{y\in L^2(\ooo);\frac{\pp y}{\pp x_i}\in L^2(\ooo),i=1,...,d\right\}$ with the norm $\|y\|_V=(|y|^2_2+|\nabla y|^2_2)^{\frac12}.$ Here, $\nabla y=\left\{\frac{\pp y}{\pp x_i}\right\}^d_{i=1}$ are taken in the sense of the  distributions space $\cald'(\ooo)$ on $\ooo$. We denote also by $V'$ the dual $(H^1(\ooo))'$ of $H^1(\ooo)$.

$W^{2\ell,\ell}_q(\calq)$ is the space of all $y\in L^q(\calq)$ such that $\frac{\pp^r}{\pp t^r}\ \frac{\pp^s}{\pp x_j}\ y\in L^q,$ $j=1,...,d$, for any $r,s$ such that $2r+s\le2\ell$. (Here $\ell$ is integer.) Similarly, there are defined the Sobolev spaces $W^{2\ell,\ell}_q(\Sigma)$ on $\Sigma$.  These definitions extend to the noninteger $\ell>0.$ (See, e.g., \cite{10}, p.~81.)

We denote by $W^{1,\9}(\calq)$ the space $\{y\in L^\9(\calq;\frac{\pp y}{\pp t},\frac{\pp y}{\pp x_i}\in L^\9(\calq),$\break $ i=1,...,d\}.$ (Here and everywhere in the following,  the derivatives are taken in the sense of distributions.)

Given a Banach space $X$ with the norm $\|\cdot\|_X$, we denote by $C([0,T];X)$ the space of $X$-valued continuous functions $y:[0,T]\to X$ with the standard supremum norm  $\|\cdot\|_{C([0,T];X)}$. By $L^2(0,T;X)$, we denote the space of \mbox{$X$-valued} measurable functions $y:(0,T)\to X$ such that $\|y\|_X\in L^2(0,T).$ By $W^{1,2}([0,T];V')$ we denote the infinite dimensional Sobolev space $$\left\{y\in L^2(0,T;V'),\frac{dy}{dt}\in L^2(0,T;V')\right\},$$ where $\frac d{dt}$ is taken in the sense of $V'$-valued distributions on $(0,T)$.

\section{The Stefan two-phase system with mushy region}\label{s2}
\setcounter{equation}{0}

We shall study here a two-phase transition model of the form \eqref{e1.1} with mushy region. Namely, we can consider the system
\begin{equation}
\label{e2.1}
\barr{ll}
\dd\frac{\pp y}{\pp t}-{\rm div}(H_\lbb(y)\nabla y)=0&\mbox{ in }\calq,\vsp
H_\lbb(y)\nabla y\cdot\nu=u&\mbox{ on }\Sigma,\vsp
y(0)=y_0&\mbox{ in }\ooo,\earr\end{equation}where
\begin{equation}
\label{e2.2}
H_\lbb(y)(t,x)=\frac1{\lbb}
\int^{t+\lbb}_tds\int_\ooo h_\lbb(y(s,x-\lbb\xi))
\vf(\xi)d\xi,\
(t,x)\in\calq.\end{equation}Here, $\lbb>0$, $\vf\in C^\9_0(\rr^d)$ is a standard mollifier, that is,
$$\vf(x)=0\mbox{ for }|x|\ge1,\ \ \vf\ge0,\ \ \int_{\rr^d}\vf(x)dx=1,$$and $h_\lbb:\rr\to\rr$ is a smooth approximation of $\beta'$. More precisely,
\begin{equation}
\label{e2.3}
h_\lbb(r)=\left\{\barr{lll}
k_1&\mbox{ for }r<-\lbb^{\alpha},\vsp
	k_2&\mbox{ for }r>\rho+\lbb^{\alpha},\vsp g_\lbb(r)&\mbox{ for }-\lbb^{\alpha}\le r\le \rho+\lbb^{\alpha},\earr\right.\end{equation}where $g_\lbb\in C^1[-\lbb^{\alpha},\rho+\lbb^{\alpha}]$ is such that
\begin{eqnarray}
&g_\lbb(-\lbb^{\alpha})=k_1,\ g_\lbb(\rho+\lbb^{\alpha})=k_2,\ g_\lbb(r)=\lbb^{\alpha},\ \ff\,r\in[0,\rho],\label{e2.4}\\[2mm]
&|g'_\lbb(r)|\le k^*\lbb^{-\alpha},\ \ \lbb^{\alpha}\le g_\lbb(r)\le k^*,\ \ff\,r\in[-\lbb^{\alpha},\rho+\lbb^{\alpha}].\label{e2.5}
\end{eqnarray}
Here $k^*=\max\{k_1,k_2\}$ and $\alpha\in(0,1)$ is a parameter which will be made precise later on. Every\-where in the following, $\lbb$ is positive and sufficiently small such that \mbox{$k^*-\lbb^{\alpha}>0$.}

By \eqref{e2.3}-\eqref{e2.5} we see that, for $\lbb\to0$,
$$H_\lbb(y)\to\beta'(y)\mbox{ in }L^1(\calq),\ \ \ff\,y\in L^2(\calq),$$and so, formally, we may view \eqref{e2.1} as an approximation to equation \eqref{e1.5}. However, the problem of the convergence of solutions to \eqref{e2.1} to the solution $y$ to \eqref{e1.5} though is interesting in itself will not be addressed here. We shall see later in (Proposition \ref{p2.3}) that equation \eqref{e2.1} is a model of the Stefan problem with mushy region discussed in Section \ref{s1}.

We shall discuss first the existence for equation \eqref{e2.1} and its relationship with the two-phase Stefan problem.

\begin{definition}\label{d2.1} \rm The function $y:\calq\to\rr$ is said to be a weak solution to equation \eqref{e2.1} if
\begin{equation}
\label{e2.6}
y\in C([0,T];L^2(\ooo))\cap L^2(0,T;V)\cap W^{1,2}([0,T];V'),\end{equation}
\begin{equation}
\label{e2.7}
\barr{r}
\dd\frac d{dt}\int_\ooo y(t,x)\psi(x)dx+\int_\ooo H_\lbb(y)(t,x)\nabla y(t,x)\cdot\nabla\psi(x)dx\\
=\dd\int_\Gamma u(t,x)\psi)(x)dx,\ \ff\psi\in V,\ \mbox{ i.e. }t\in(0,T),\earr\end{equation}
\begin{equation}
\label{e2.8}
y(0,x)=y_0(x),\mbox{\ \ a.e. }x\in\ooo,\end{equation}
where $V=H^1(\ooo),\ V'=(H^1(\ooo))'.$
\end{definition}

\begin{proposition}\label{p2.2} Let $y_0\in L^2(\ooo)$, $\lbb>0$, and $u\in y=y^u$. Then there is a unique weak solution $y=y^u$ to equation \eqref{e2.1} and one has
\begin{equation}
\label{e2.9}
|y(t)|^2_2+\lbb^{\alpha}\int^t_0\|y(s)\|^2_Vds\le C(|y_0|^2_2+\lbb^{-\alpha}\|u\|^2_{L^2(\sigma)}),\end{equation}
\begin{equation}
\label{e2.9a}
\int^T_0\left\|\frac{dy}{dt}\,(t)\right\|^2_{V'}dt
\le C\lbb^{-\alpha}(|y_0|^2_2+\lbb^{-\alpha}
\|u\|^2_{L^2(\Sigma)}),\end{equation}where $C$ is independent of $\lbb$. Assume further that $y_0\in V$. Then
\begin{equation}
\label{e2.10}
 \barr{r}
\dd
|y(t,x){-}y(s,\xi)|\le C(\ooo_0)\lbb^{-\frac{13\alpha}2}
(\|y_0\|_V{+}\|u\|_{L^2(\Sigma)})
(|t{-}s|^{\frac12}{+}|x{-}\xi|^{\frac12}),\vsp
\ff(t,x),(s,\xi)\in\calq_0=[0,T]\times\ooo_0,\earr \end{equation}where $\ooo_0$ is any open set such that $\overline{\ooo}_0\subset\ooo_0.$
\end{proposition}

\pf $1^\circ$. {\it Existence.} We set
$$\barr{l}
K_M=\Bigg\{ z\in C([0,T];L^2(\ooo))\cap L^2(0,T;V)\cap W^{1,2}([0,T];V');\\
\hspace*{18mm}\|z\|_{C([0,T];L^2(\ooo))}+\|z\|_{L^2(0,T;V)}
+\left\|\dd\frac{dz}{dt}\right\|_{L^2(0,T;V')}\le M\Bigg\}.\earr$$
For each $z\in K_M$, the linear parabolic equation
\begin{equation}
\label{e2.11}
\barr{ll}
\dd\frac{\pp y}{\pp t}-{\rm div}(H_\lbb(z)\nabla y)=0&\mbox{ in }\calq,\vsp
H_\lbb(z)\nabla y,\nu=u&\mbox{ on }\Sigma,\vsp
y(0)=y_0&\mbox{ in }\ooo,\earr
\end{equation}has a unique solution $y=\Phi(z)\in C([0,T];L^2(\ooo))\cap L^2(0,T;V)\cap W^{1,2}([0,T];V'),$ that is,

\begin{equation}
\label{e2.12}
\barr{r}
\dd\frac d{dt}\int_\ooo y(t,x)\psi(x)dx+\dd\int_\ooo H_\lbb(z)(t,x)\nabla y(t,x)\cdot\nabla\psi(x)dx\\
=\dd\int_\Gamma u(t,x)\psi(x)dx,\ \ \ff\psi\in V,\ \mbox{a.e. }t\in(0,T).\earr\end{equation}Moreover, recalling that $H_\lbb\ge \lbb^{\alpha}$, a.e. on $\calq$, we get by \eqref{e2.12}$$\frac12\,|y(t)|^2_2+\lbb^{\alpha}\int^t_0|\nabla y(s)|^2_2ds
\le\frac12\,|y_0|^2_2+\int^t_0\int_\Gamma u(s,x)y(s,z)ds\,dx$$and, by the trace theorem, this yields
$$\barr{l}
\|y(t)\|^2_{C([0,T];L^2(\ooo))}+2\lbb^{\alpha}\|\nabla y|^2_{L^2(0,T;L^2(\ooo))}\vsp\qquad
\le|y_0|^2_2+\lbb^{-\alpha}\|u\|^2_{L^2(\Sigma)}+\lbb^{\alpha}\|y\|^2_{L^2(0,t;L^2(\Gamma))},\ \ff\,t\in[0,T],\earr$$and so, by the trace theorem we get
\begin{equation}
\label{e2.13}
\|y(t)\|^2_{C([0,T];L^2(\ooo))}+\lbb^{\alpha}\int^t_0\|y(s)\|^2_V\le C(|y_0|^2_2+\lbb^{-\alpha}\|u\|^2_{L^2(\Sigma)}).\end{equation}
Moreover, since $\|H_\lbb(t)\|_{L^\9(\calq)}\le k^*$, we see by \eqref{e2.12} that
$$\left\|\frac d{dt}\,y(t)\right\|_{V'}\le C(|\nabla y(t)|_2+\|u(t)\|_{L^2(\Gamma)}),\ \ff\,t\in[0,T],$$and, therefore, by \eqref{e2.13} we have
\begin{equation}\label{e2.13a}\int^T_0\left\|\frac d{dt}\,y\right\|^2_{V'}dt\le C\lbb^{-\alpha}(|y_0|^2_2+\lbb^{-\alpha}\|u\|^2_{L^2(\Sigma)}),\end{equation}where $C$ is independent of $\lbb$.
Hence, for $M$ sufficiently large, $\Phi(K_M)\subset K_M$. Since $\phi$ is continuous in $L^2(\calq)$-norm
and by the Aubin-Lions compactness theorem (see, e.g., \cite{11}, p.~61) $K_M$ is compact in $^2(\calq)$, by Schauder's fixed point theorem there is $y\in K_M$ such that $\Phi(y)=y.$ Clearly, $y$ is a solution to \eqref{e2.1} and also estimate \eqref{e2.9} follows.

\bk\n$2^\circ$. {\it Uniqueness.} If $y_1,y_1$ are two weak solutions to \eqref{e2.1}, we have, for $y=y_1-y_2$,

$$\barr{ll}
\dd\frac {\pp y}{\pp t}-{\rm div}(H_\lbb(y_1)\nabla y)-{\rm div}(H_\lbb(y)\nabla y_1)=0&\mbox{ in }\calq,\vsp
(H_\lbb(y_1)\nabla y+H_\lbb(y)\nabla y_1)\cdot\nu=0&\mbox{ on }\Sigma,\vsp
y(0)=0&\mbox{ in }\ooo.\earr$$
This yields
$$\barr{l}
\dd\frac 12\ \frac d{dt}\ |y(t)|^2_2+\int_\ooo H_\lbb(y_1)(t,x)|\nabla y(t,x)|^2dx\vsp
\qquad=\dd\int_\ooo H_\lbb(y)(t,x)\nabla y_1(t,x)\cdot\nabla y(t,x)dx,\ \mbox{ a.e. }t\in(0,T).\earr$$
Since $H_\lbb\ge\lbb^{-\alpha},$ a.e. in $\calq$, we get
$$|y(t)|^2_2+2\lbb^{-\alpha}\int^t_0|\nabla y(s)|^2_2dx\le C\int^t_0|y(s)|_2|\nabla y(s)|_2ds,\ \ff \,t\in(0,T),$$
which implies $y\equiv0$, as claimed.

\bk
\n$3^\circ.$ {\it Regularity.} Assume now that $y_0\in H^1(\ooo)$ and consider the arbitrary open subsets $\ooo_1,\ooo_2$ of $\ooo$ such that $\ov\ooo_0\subset\ooo_1$, $\ov\ooo_1\subset\ooo_2,$ $\ov\ooo_2\subset\ooo.$ Assume that the boundary $\pp\ooo_0$ as well as that of $\ooo_1,\ooo_2$ are  smooth (of class $C^2$, for instance).

Let $\calx\in C^\9_0(\ooo)$ be such that $\calx=1$ on $\ov\ooo_2$ and let $\wt y=\calx y$. We have by \eqref{e2.1} that
\begin{equation}
\label{e2.14}
\barr{ll}
\dd\frac{\pp\wt y}{\pp t}-H_\lbb(y)\Delta\wt y=g&\mbox{ in }\calq,\vsp
\wt y=0&\mbox{ on }\Sigma,\vsp
\wt y(0)=y_0\calx&\mbox{ in }\ooo,\earr\end{equation}
where
\begin{equation}
\label{e2.15}g=-y\nabla H_\lbb(y)\cdot\nabla\calx-2H_\lbb(y)\nabla y\cdot\nabla\calx-H_\lbb(y)y\Delta\calx.\end{equation}To this end, we shall use a bootstrap argument and that \eqref{e2.10} holds.
We are going to prove that $y$ is H\"older continuous of order $\frac12$ on $\calq_0=(0,T)\times\ooo_0$. Namely, by \eqref{e2.14} and \eqref{e2.15}, we shall prove first that
\begin{equation}
\label{e2.18}
\|y\|_{W^{2,1}_{\frac32}}(\calq)\le C(\ooo_0,y_0,u)\lbb^{-2\alpha},\ \ff\lbb\in(0,k^*)\end{equation}and we shall use this estimate to show that
\begin{equation}
\label{e2.19}
\|y\|_{W^{2,1}_2(\calq_1)}\le C(\ooo_0,y_0,u)\lbb^{-6\alpha}\end{equation}and, finally, that\newpage 
\begin{equation}
\label{e2.19a}
\|y\|_{W^{2,1}_3(\calq_0)}\le C(\ooo_0,y_0,u)\lbb^{-\frac{13\alpha}2},\end{equation} which will lead to the desired estimate \eqref{e2.10}. Here and everywhere in the following,
$$C(\ooo_0,y_0,u)=CC^*(\|y_0\|_V+\|u_0\|_{L^2(\Sigma)}),$$where $C^*=C(\ooo_0,\ooo_1,\ooo_2)$ is fixed (but depends of $\ooo_0,\ooo_1,\ooo_2$) and $C$ is a constant which is independent of $\lbb,$ $\ooo_i$, $i=0,1,2,$ but will change from one estimate to another.

We note first that, by \eqref{e2.5} and \eqref{e2.15}, we have
\begin{equation}
\label{e2.20}
|g|\le C^*(|y|\,|\nabla H_\lbb(y)|+|\nabla y|+|y|),\mbox{ a.e. in }\calq.\end{equation}
On the  other hand, by \eqref{e2.2} and \eqref{e2.5} we see that
\begin{equation}
\label{e2.21}
\barr{ll}
|\nabla H_\lbb(y)(t,x)|\!\!\!
&C\lbb^{-1-\alpha}
\dd\int^{t+\lbb}_tds\int_{[|\xi|\le1]}|\nabla y(s,x-\lbb\xi)|d\xi\vsp
&\le C\lbb^{-1-d-\alpha}
\dd\int^{t+\lbb}_tds\int_{[|x-\xi|\le\lbb]}|\nabla y(s,\xi)|d\xi\vsp
&\le C\lbb^{-\frac d2-\alpha}\(\dd\frac1\lbb
\dd\int^{t+\lbb}_tds\int_{[|x-\xi|\le\lbb]}|\nabla y(s,\xi)|^2d\xi\)^{\frac12}.
\earr\end{equation}
This yields
\begin{equation}
\barr{l}
\dd\int_\calq|\nabla H_\lbb(y)(t,x)|^2dt\,dx\\
\qquad\qquad\le C\lbb^{-2\alpha}
\dd\int^T_0\frac1\lbb
\dd\int^{t+\lbb}_tds\int_{[|\xi-\xi|\le\lbb]}|\nabla y(s,\xi)|^2d\xi\,dt\vsp
\qquad\qquad\le C\lbb^{-2\alpha}\dd\int_\calq|\nabla y(t,x)|^2dt\,dx.\earr
\label{e2.22}
\end{equation}
On the other hand, by \eqref{e2.9} we have
\begin{equation}
\label{e2.23}
|y(t)|^2_2+\dd\int^T_0\|y(t)\|^2_Vdt
\le C\lbb^{-2\alpha}(|y_0|^2_2+\|u\|^2_{L^2(\Sigma)})
\le(C(\ooo_0,y_0,u))^2\lbb^{-2\alpha}\end{equation}
and so, by \eqref{e2.22}, this yields
\begin{equation}
\label{e2.24}
\|\nabla H_\lbb(y)\|_{L^2(\calq)}\le C(\ooo_0,y_0,u)\lbb^{-2\alpha}.\end{equation}
Then, by \eqref{e2.20}, \eqref{e2.22}, \eqref{e2.23}, we get by the Sobolev embedding theorem combined with H\"older's inequality
$$\barr{ll}
\dd\int_\calq|g|^{\frac32}dt\,dx\!\!\!
&\le\dd C\int_\calq|\nabla H_\lbb(y)|^{\frac32}|y|^{\frac32}dt\,dx+C\dd\int_\calq(|\nabla y|^{\frac32}+|y|^{\frac32})dt\,dx\vsp
&\le \dd\int^T_0|\nabla H_\lbb(t)|^{\frac32}_2|y(t)|^{\frac32}_6dt
+(C(\ooo_),y_0,u))^{\frac32}\lbb^{-\frac{3\alpha}2}\vsp
&\le(C(\ooo_0,y_0,u))^{\frac32}\lbb^{-\frac{3\alpha}2}+\|\nabla H_\lbb(y)\|^{\frac32}_{L^2(\calq)}
\|y(t)\|^{\frac32}_V\vsp
&\le (C(\ooo_0,y_0,u))^{\frac32}\lbb^{-3\alpha}.
\earr$$Hence, by \eqref{e2.20}, it follows that $g\in L^{\frac32}(\calq)$ and
$$\|g\|_{L^{\frac23}}(\calq)\le C(\ooo_0,y_0,u)^{-2\alpha}.$$
	Then, by \eqref{e2.14} it follows that $\wt y\in W^{2,1}_{\frac32}(\calq)$ and (see, e.g., \cite{10}, p.~342)
	$$\|\wt y\|_{W^{1,2}_{\frac32}(\calq)}\le C\lbb^{-\alpha}(\|g\|_{L^{\frac32}(\calq)}+\|\calx y_0\|_{W^{\frac23}_{\frac32}(\ooo)})\le C(\ooo_0,y_0,u)\lbb^{-3\alpha}.$$
	Taking into account that $\wt y=y$ on $\calq_2$, we get \eqref{e2.18}, as claimed.
	
	For simplicity, in the following we shall just write $C$ instead of $C(\ooo_0,y_0,u)$. To prove \eqref{e2.19}, we note first that, since $W^{2,1}_{\frac32}(\calq_2)\subset L^\9(\calq_2)$ in $3-D$, we have by \eqref{e2.18}
\begin{equation}
\label{e2.25}
\|y\|_{L^\9(\calq_2)}\le C\lbb^{-3\alpha}.\end{equation}Now, we denote again by $\wt y$ the function $\calx y$,  where $\calx\in C^\9_0(\ooo_2)$ is such that $\calx=1$ on $\ov\ooo_1.$ Taking into account \eqref{e2.20} and \eqref{e2.24}, we get
$$|g|\le C(|\nabla H_\lbb(y)|+1)\lbb^{-3\alpha}+|\nabla y|,\mbox{ a.e. in }\calq_2,$$ and, therefore,  by \eqref{e2.23}
$$\|g\|_{L^2(\calq_2)}\le C\lbb^{-3\alpha} (\|\nabla H_\lbb(y)\|_{L^2(\calq_2)}+1)+C\lbb^{-\alpha}.$$On the other hand, by \eqref{e2.21} we get, as above (see \eqref{e2.22}-\eqref{e2.24}),
$$\|\nabla H_\lbb(y)\|_{L^2(\calq_2)}\le C\lbb^{-2\alpha}.$$This yields
$$\|g\|_{L^2(\calq_2)}\le C\lbb^{-5\alpha}$$and, again by \eqref{e2.14}, we infer that $\wt y=W^{2,1}_2(\calq_2)$ and
$$\|\wt y\|_{W^{2,1}_2(\calq_2)}\le C\lbb^{-6\alpha}$$and so \eqref{e2.19} follows.

We also note for later use that, by \eqref{e2.14} taken in $\calq_2$, it also follows
$$\|\nabla\wt y(t)\|_{L^2(\ooo_2)}+\lbb^{-\alpha}\int^t_0\|\Delta\wt y(s)\|^2_{L^2(\ooo_2)}ds
\le\|g\|_{L^2(\calq_2)}+\|y_0\calx\|_{W^{1}_2(\ooo_2)}$$and so
\begin{equation}
\label{e2.26}
\|\nabla y(t)\|_{L^2(\ooo_1)}\le\|\nabla \wt y(t)\|_{L^2(\ooo_2)}\le C\lbb^{-5\alpha},\ \ \ff\,t\in[0,t].\end{equation}

Now, we take $\wt y=\calx y$, where $\calx\in C^\9_0(\ooo_1)$ and $\calx=1$ on $\ooo_0.$ By \eqref{e2.25}, we have as above
$$|g|\le C(|\nabla H_\lbb(y)|+1)\lbb^{-3\alpha}+|\nabla y|,\mbox{ a.e. in }\calq_1,$$ and, therefore,
\begin{equation}
\label{e2.27}
\|g\|_{L^3(\calq_1)}\le C\lbb^{-3\alpha}(\|\nabla H_\lbb(y)\|_{L^3(\calq_1)}+1)+C\|\nabla y\|_{L^3(\calq_1)}.\end{equation}
On the other hand, by interpolating between $L^2$ and $L^6$, we have
$$\|\nabla y(t)\|_{L^3(\ooo_1)}\le\|\nabla y(t)\|^{\frac34}_{L^2(\ooo_1)}
\|\nabla y(t)\|^{\frac14}_{L^6(\ooo_1)}
\le C\|\nabla y(t)\|^{\frac34}_{L^2(\ooo_1)}
\|\nabla y(t)\|^{\frac14}_{H^1(\ooo_1)}.$$
By \eqref{e2.23}, \eqref{e2.26}, this yields
$$\dd\int^T_0\|\nabla y(t)\|^3_{L^3(\ooo_1)}dt
\le T\|\nabla y\|^4_{L^\9(0,T;L^2(\ooo_1))}
\dd\int^T_0\|y(t)\|^{\frac34}_{H_1(\ooo_1)}dt
\le C\lbb^{-\frac{9\alpha}2}.$$Equivalently,
\begin{equation}
\label{e2.28}
\|\nabla y\|_{L^3(\calq_1)}\le C\lbb^{-\frac{3\alpha}2}.
\end{equation}
On the other hand, by \eqref{e2.21} we have
$$\barr{ll}
|\nabla H_\lbb(y)(t,x)|^3\!\!\!
&\le C\lbb^{-3(\alpha+d+1)}
\(\dd\int^{t+lbb}_tds\int_{[|x-\xi|\le\lbb]}|\nabla y(s,\xi)|d\xi\)^3\vsp
&\le C\lbb^{-3(\alpha+1)-d}\(\dd\frac1\lbb\int^{y+\lbb}_tds\(\int_{[|x-\xi|\le\lbb]}|\nabla y(s,\xi)|^3d\xi\)^{\frac13}\)^3\vsp
&\le C\lbb^{-3\alpha-d}\dd\frac1\lbb\int^{t+\lbb}_tds\int_{[|x-0\xi|\le\lbb]}|\nabla y(s,\xi)|^3d\xi,\ \ff(t,x)\in\calq_1.\earr$$
This yields (see \eqref{e2.28})
$$\|\nabla H_\lbb(y)\|_{L^3(\calq_1)}\le C\lbb^{-\alpha}\|\nabla y\|_{L^3(\calq_1)}\le C\lbb^{-\frac{5\alpha}2}$$and so, by \eqref{e2.27}, we have
$$\|g\|_{L^3(\calq_1)}\le C\lbb^{-\frac{11\alpha}2}.$$Then, once again by equation \eqref{e2.14} taken in $\calq_1$, we infer that $\wt y\in  W^{2,1}_3(\calq_1)$~and (see \cite{10}, p.~342)
$$\|\wt y\|_{W^{2,1}_3(\calq_1)}\le C\lbb^{-\frac{13\alpha}2},$$which implies \eqref{e2.19a}, as claimed.

Then, by \eqref{e2.19a}  and Lemma 3.3 in \cite{10}, p.~80, we infer  that $y$ is H\"older continuous of order $\frac12$ on $\calq_0$ and  that \eqref{e2.10} follows.\hfill$\Box$\bk

Let  $\ooo_0$ be as above an open subdomain of $\ooo$ with smooth boundary $\pp\ooo_0$ such that $\ov\ooo_0\subset\ooo$ and $\lbb\in(0,1)$ sufficiently small such that 
$${\rm dist}(\ov\ooo_0,\Gamma)\ge\lbb.$$
We have

\begin{proposition}\label{p2.3} Let $y=y^u$ be the solution to equation \eqref{e2.1} for $u\in L^2(\Sigma)$ and $y_0\in H^1(\ooo)$.  Then, for  $\alpha\in(0,(26)^{-1})$, we have	
\begin{eqnarray}
\label{e2.29}
\dd\frac{\pp y}{\pp t}-k_1\Delta y=0&\mbox{ in }&\{(t,x)\in\calq_0;y(t,x)\le-2\lbb^{\frac14}\},\\
\dd\frac{\pp y}{\pp t}-k_2\Delta y=0&\mbox{ in }&\{(t,x)\in\calq_0;y(t,x)\ge\rho+2\lbb^{\frac14}\}.\label{e2.30}
\end{eqnarray}\end{proposition}
This means that $\sigma^u_t=\{x\in\ooo;\ -2\lbb^{\frac14}\le y^u(t,x)\le\rho+2\lbb^{\frac14}\}$ can be viewed as a mushy set of system \eqref{e1.1}.\bk

\n{\bf Proof.} By \eqref{e2.2}, we have
$$\barr{ll}
H_\lbb(y)(t,x)\!\!\!
&=\dd\frac1\lbb\int^{t+\lbb}_tds\int_{[|\xi|\le1]}h_\lbb(y(s,x-\lbb\xi))\vf(\xi)d\xi\vsp
&=h_\lbb(y(t,x))+F_\lbb(t,x),\ \ff(t,x)\in\calq_0,\earr$$where
$$|F_\lbb(t,x)|\le\frac1\lbb\int_{[|t-s|\le\lbb;|\xi|\le1]}
|h_\lbb(s,x-\lbb\xi)-h_\lbb(y(t,x))|ds\,d\xi.$$By \eqref{e2.10}, we have
$$|y(s,x-\lbb\xi)-y(t,x)|\le C(\ooo_0)(|t-s|^{\frac12}+|x-\xi|^{\frac12})\lbb^{-\frac{13\alpha}2}$$for $(x,t)\in\calq_0,\ |\xi|\le1.$ Taking into account that $13\alpha<\frac12$,   this yields
$F_\lbb(t,x)=0$ and $H_\lbb(y)(t,x)=h_\lbb(y(t,x))=k_1$ for $y(t,x)\le-2\lbb^{\frac14}$ and $\lbb$ sufficiently small. Similarly, follows \eqref{e2.30}.\hfill$\Box$\bk

In terms of the temperature $\theta$ (as in \eqref{e1.1}), it follows by \eqref{e2.18}-\eqref{e2.19} that
$$\barr{ll}
\dd\frac{\pp\theta}{\pp t}-k_1\Delta\theta=0&\mbox{ in }\{(t,x)\in \calq_0;\theta(t,x)\le-2\lbb^{\frac14}\}\vsp
\dd\frac{\pp\theta}{\pp t}-k_1\Delta\theta=0&\mbox{ in }\{(t,x)\in\calq_0;\theta(t,x)\ge\rho+2\lbb^{\frac14}\},\earr$$while $\nabla\theta=H_\lbb(y)\nabla y$ in
\begin{equation*}
\sigma^u_t=\{x\in\ooo_0;|\theta(t,x)|\le2\lbb^{\frac14}\}\cup(\ooo\setminus\ooo_0),\end{equation*}where $y=y^u$ is the solution to \eqref{e2.1}. This means that $\sigma^u_t$ can be viewed as the {\it mushy region} of the phase transition process described by equation \eqref{e2.1}. It should be emphasized that in formulation \eqref{e1.5} of equation \eqref{e1.1} the fact that the function $\beta^{-1}$  is discontinuous in origin implies that the flux $q$ of $\theta$ has a jump $\rho$ on the free boundary $S=\{\theta=0\}$. By replacing  in \eqref{e1.1} the function $\beta(y)$ by a smooth version $H_\lbb(y)$, the temperature flux $q$ is replaced by a "mild" phase transition with a mushy region $\sigma^u_t$ which separates the solid and liquid regions.

One also should mention   that \eqref{e2.1} is only one  of possible models of the phase transition with mushy region. Another popular model  is that intro\-duced by G. Caginalp \cite{5}, which is described by the phase-field system and from which one gets to the limit problem \eqref{e1.1}. The model considered here, that is \eqref{e2.1}, is different and remains on the physical grounds of two-phase transition mechanism.

\section{The controllability result}\label{s3}
\setcounter{equation}{0}

Everywhere in the following, the constants $\lbb$ and $\mu$ are positive and arbitra\-rily small but fixed. For $u\in L^2(\Sigma)$, denote as above by $y^u$ the solution to a system corresponding to $y_0\in L^2(\ooo)$. Finally, $\sigma^u_t$ is the  set
\begin{equation}
\label{e3.1}
\sigma^u_t=\{x\in\ooo;-2\mu\le y^u(t,x)\le\rho+2\mu\},\ t\in[0,T].\end{equation}

\begin{theorem}\label{t3.1} Let $T>0$ and $\ooo^*$ be a Lebesgue measurable subset of $\ooo$ with positive measure. Then, under assumptions \eqref{e2.3}-\eqref{e2.5}, there is a controller $u^*\in L^2(\Sigma)$ such that
\begin{equation}
\label{e3.2}
\ooo^*\subset\sigma^{u^*}_T.\end{equation}
\end{theorem}

We note that, for $\mu=2\lbb^{\frac14}$ and  $0<\alpha<(26)^{-1}$, it follows by Proposition \ref{p2.3} that $\sigma^u_t$ is just the mushy set $\{x\subset\ooo;\ -2\lbb^{\frac14}\le y^u(t,x)\le\rho+2\lbb^{\frac14}\}$ of system \eqref{e1.1}. Then, by Theorem \ref{t3.1}, we have

\begin{corollary}\label{c3.2} Assume that $\alpha\in(0,(26)^{-1}).$ Then, for every Lebesgue measurable set $\ooo^*\subset\ooo$ with positive measure, there is $u^*\in L^2(\Sigma)$ such that
	$$\ooo^*\subset\{x\in\ooo;\ -2\lbb^{\frac14}\le y^{u^*}(T,x)\le\rho+2\lbb^{\frac14}\}.$$\end{corollary}

\n{\bf Proof of Theorem \ref{t3.1}.} We shall prove first \eqref{e3.2} for the linear control system
\begin{equation}
\label{e3.3}
\barr{ll}
\dd\frac{\pp y}{\pp t}-{\rm div}(H_\lbb(z)\nabla y)=0&\mbox{ in }\calq,\vsp
H_\lbb(z)\nabla y\cdot\nu=u&\mbox{ on }\Sigma,\vsp
y(0)=y_0&\mbox{ in }\ooo,\earr\end{equation}where $z\in L^2(0,T;V)$ is arbitrary but fixed.

Denote by $y^{u,z}\in L^2(0,T;V)\cap C([0,T];L^2(\ooo);V)\cap W^{1,2}([0,T];V')$ the solution to \eqref{e3.3}. We have

\begin{lemma}\label{l3.2} There is $u^*\in L^2(\Sigma)$ such that
\begin{equation}
\label{e3.4}
\ooo^*\subset\sigma^{u^*,z}_T=\{x\in\ooo;-\mu\le y^{u^*,z}(T,x)\le\rho+\mu\}.\end{equation}
Moreover, one has
\begin{equation}
\label{e3.5}
|y^{u^*,z}(t)|^2_2+\lbb^{\alpha}\dd\int^t_0|\nabla y^{u^*,z}(s)|^2_2ds
\le C\(|y_0|^2_2+\lbb^{-\alpha}\dd\int^t_0\int_\Gamma|u^*|^2ds\,dx\),\end{equation}
\begin{equation}
\label{e3.6}
\dd\int^T_0\left|\frac d{dt}\ y^{u^*,z}(t)\right|^2_{V'}dt
\le C\lbb^{-\alpha}(|y_0|^2_2
+\lbb^{-\alpha}\|u\|^2_{L^2(\Sigma)}),\end{equation}
\begin{equation}
\label{e3.7}
\|u^*\|_{L^2(\Sigma)}\le C\lbb^{-\alpha}|u_0|_2,\end{equation}
where $C$ is independent of $\lbb$.
	\end{lemma}

\n{\bf Proof.} For each $\vp>0$ consider the minimization problem
\begin{equation}
\label{e3.8}
\barr{r}
{\rm Min}\Big\{\dd\frac12\ \|u\|^2_{L^2(\Sigma)}+\dd\frac1{2\vp}\dd\int_{\ooo^*}((y^{u,z}(T,x)+\mu)^-)^2\vsp
+((y^{u,z}(T,x)-\rho-\mu)^+)^2dx\Big\}.\earr\end{equation}
It is easily seen by estimates \eqref{e2.9}, \eqref{e2.9a} that there is a unique solution $u_\vp$ to problem \eqref{e3.8}. We set $y^{u,z}=y_\vp$ and note that by the first order optimality conditions in \eqref{e3.8} we have
\begin{equation}
\label{e3.9}
u_\vp=p_\vp\mbox{ on }\Sigma,\end{equation}where $p_\vp\in L^2(0,T;V)\cap C([0,T];L^2(\ooo))\cap W^{1,2}([0,T];V')$ is the solution to the dual backward system
\begin{equation}
\label{e3.10}
\barr{ll}
\dd\frac{\pp p_\vp}{\pp t}+{\rm div}(H_\lbb(z)\nabla p_\vp)=0&\mbox{ in }\calq,\vsp
\dd\frac{\pp p_\vp}{\pp\nu}=0&\mbox{ on }\Sigma,\vsp
p_\vp(t)=\dd\frac1\vp\ {\bf1}_{\ooo^*}
((y_\vp(T)+\mu)^--(y_\vp(T)-\rho-\mu)^+)&\mbox{ in }\ooo,\earr \end{equation}
where ${\bf1}_{\ooo^*}$ is the characteristic function of $\ooo^*$. By \eqref{e3.3}, where $u=u_\vp,$ $ y^{u,z}=y_\vp$, and by system \eqref{e3.10} we get
$$\barr{ll}
\!\!\dd\int_\Sigma u^2_\vp\,dt\,dx=\!\!\!
&-\!\!\dd\int_\ooo p_\vp(0,x)y_0(x)dx+\!\!\int_\ooo p_\vp(T,x)y_\vp(T,x)dx=-\!\!\dd\int_\ooo p_\vp(0,x)y_0(x)\vsp&+\dd\frac1\vp\int_{\ooo^*}y_\vp(T,x)((y_\vp(T,x)+\mu)^--(y_\vp(T,x)-\rho-\mu)^+)dx.\earr$$This yields
\begin{equation}
\label{e3.11}
\hspace*{-5mm}\barr{r}
\dd\int_\Sigma u^2_\vp\,dt\,dx+\frac1\vp\int_{\ooo^*}(((y_\vp(T,x)+\mu)^-)^2+((y_\vp(T,x)-\rho-\mu)^+)^2dx\\
\le\dd\int_\ooo p_\vp(0,x)y_0(x)dx
\le|p_\vp(0)|_2|y_0|_2,\ \ \ff\,\vp>0.\earr\hspace*{-2mm}\end{equation}On the other hand, we have the following Carleman's type inequality for the solution $p_\vp$ to equation \eqref{e3.10}
\begin{equation}
\label{e3.12}
|p_\vp(0)|_2\le C\lbb^{-\alpha}\|p_\vp\|_{L^2(\Sigma)},\ \ \ff\,\vp>0,\end{equation}where $C$ is independent of $\vp$ and $\lbb$.

Here is the proof. It is well known that, for each $\bar y_0\in L^2(\ooo)$, there is a boundary controller $u\in L^2(\Sigma)$ such that the corresponding solution $\bar y^{u,z}$ in $T$ and
\begin{equation}
\label{e3.13}
\|u\|_{L^2(\Sigma)}\le C\|H_\lbb(z)\|_{W^{1,\9}(\calq)}|y_0|_2\le C\lbb^{-1}|\bar y_0|_2,\ \ff\bar y_0\in L^2(\ooo),\end{equation}where $C$ is independent of $z$ and $y_0$ because,
  by \eqref{e2.2} and \eqref{e2.5}, we have the estimate
$$\|H_\lbb(t)\|_{L^\9(\calq)}+\|\nabla H_\lbb(t)\|_{L^\9(\calq)}\le k^*(m(\ooo)+\lbb^{-\alpha}\|\nabla z\|_{L^2(\calq)}).$$
In~fact, the above boundary controllability result follows in a standard way by the Carleman inequa\-lity for  linear parabolic equations with $W^{1,\9}(\calq)$ coefficients (see, e.g., \cite{2a}, \cite{7}) and applies neatly to system \eqref{e3.10}.
Then \eqref{e3.12} follows by \eqref{e3.13} by a duality argument. Indeed,  we see by \eqref{e3.10} and \eqref{e3.3} that
\begin{equation}
\label{e3.14}
\barr{ll}
\dd\int_\ooo \bar y_0(x)p_\vp(0,x)dx\!\!\!
&=\dd\int_\ooo \bar y^{u,z}(T,x)p_\vp(T,x)dx+\dd\int_\Sigma u(t,x)p_\vp(t,x)dx\vsp&=\dd\int_\Sigma u(t,x)p_\vp(t,x)dt\,dx.\earr\end{equation}Then, since $\bar y_0$  is arbitrary in $L^2(\ooo)$, it follows by \eqref{e3.13}, \eqref{e3.14} that \eqref{e3.12} holds.

Now, taking into account \eqref{e3.9} and substituting \eqref{e3.12} in \eqref{e3.11}, we get the estimate
\begin{equation}
\label{e3.15}
\hspace*{-4mm}\barr{l}
\dd\int_\Sigma|u_\vp|^2dt\,dx+\dd\frac1\vp\!\int_{\ooo^*}\!\!(((y_\vp(T,x){+}\mu)^-)^2{+}((y_\vp(T,x){-}\rho{-}\mu)^-)^2)dx\vsp\qquad\qquad\qquad
\le C\lbb^{-2\alpha}|y_0|^2_2.\earr\end{equation}Then, for $\vp\to0$, we have
\begin{equation}
\label{e3.16}
\|(y_\vp(T)+\mu)^-\|_{L^2(\ooo^*)}+\|(y_\vp(T)-\rho-\mu)^+\|_{L^2(\ooo^*)}\to0\end{equation}and, on a subsequence $\{\vp\}\to0$,
\begin{equation}
\label{e3.17}
u_\vp\to u^*\mbox{ weakly in }L^2(\Sigma).\end{equation}Moreover, by  \eqref{e2.9}, \eqref{e2.9a} we have
\begin{equation}
\label{e3.18}
\|y_\vp\|^2_{C([0,T];L^2(\ooo))}+\!\!\dd\int^t_0\!\!\(\!|\nabla y_\vp(s)|^2_2
+\left\|\dd\frac{dy_\vp}{ds}\,(s)\right\|^2_{V'}\)ds
\le C(1+\lbb^{-2\alpha})|y_0|^2_2\ \ \end{equation}
and so, again by the Aubin-Lions compactness theorem, we have on a subsequence $\{\vp\}\to0$
$$\barr{rcll}
y_\vp&\to&y^{u^*,z}&\mbox{ weakly in }L^2(0,T;V), \mbox{ strongly in }L^2(\calq)\vsp
\dd\frac{dy_\vp}{dt}&\to&\dd\frac d{dt}\,y^{u^*,z}&\mbox{ weakly in }L^2(0,T;V')\vsp
y_\vp(t)&\to&y^{u^*,z}(t)&\mbox{ weakly in }L^2(\ooo),\mbox{ uniformly in }t\in[0,T].\earr$$Then, by \eqref{e3.16}, we get
$$(y^{u^*,z}(T)+\mu)^-=0,\ \ (y^{u^*,z}(T)-\rho-\mu)^+=0,\mbox{ a.e. on }\ooo^*,$$and so \eqref{e3.4} follows. Moreover, estimates \eqref{e3.5}-\eqref{e3.7} follow by \eqref{e3.15} and \eqref{e3.18}.\hfill$\Box$\bk

 \n{\bf Proof of Theorem \ref{t3.1} (continued).} Consider the set
$$\calk=\{z\in L^2(\calq);\ \|z\|_{L^2(\calq)}\le N\},$$where $N>0$ will be made precise later on. By Lemma \ref{l3.2} we know that, for each $z\in\calk$, there is at least one $u\in L^2(\Sigma)$ such that
\begin{equation}
\label{e3.19}
-\mu\le y^{u,z}(T)\le\rho+\mu,\mbox{ a.e. in }\ooo^*,\end{equation}
\begin{equation}
\label{e3.20}
\|u\|_{L^2(\Sigma)}\le C\lbb^{-\alpha}|y_0|_2,\end{equation}where $C$ is independent of $\lbb$ and $y_0$.

We denote by $\Psi(z)\subset L^2(\calq)$ the set  of all $y^{u,z}$ satisfying \eqref{e3.19}, \eqref{e3.20}. The multivalued operator $\Psi:\calk\to L^2(\calq)$ is upper-semicontinuous. (Since $\Psi(K)$ is relatively compact, this property is equivalent with the fact that it is  strongly closed in $L^2(\calq)\times L^2(\calq)$.)    Moreover, by \eqref{e2.13}, we have
\begin{equation}
\label{e3.21}
\|y^{u,z}\|_{C([0,T];L^2(\ooo))}+\lbb^{\alpha}\int^t_0\|y^{u,z}(s)\|^2_{V}ds\le C|y_0|^2_2(1+\lbb^{-\alpha})\end{equation}and, as seen earlier (see \eqref{e2.13a}),
\begin{equation}
\label{e3.22}
\barr{ll}
\dd\int^T_0\left\|\frac d{dt}\,y^{u,z}(t)\right\|^2_{V'}dt\!\!\!&
\le C\lbb^{-\alpha}(|y_0|^2_2+\lbb^{-\alpha}\|u\|^2_{L^2(\Sigma)}\vsp&\le C\lbb^{-\alpha}(1+\lbb^{-\alpha})|y_0|^2_2.\earr\end{equation}By \eqref{e3.21}, \eqref{e3.22}, it follows that $\Psi(\calk)\subset\calk$ for $N$ sufficiently
large and also, by the compactness theorem mentioned above, that $\Psi(\calk)$ is a relatively compact set in $L^2(\calq).$ Then, by the Kakutani fixed point theorem (see, e.g., \cite{3}, p.~7) there is $y^*\in\calk$ such that $y^*\in\Psi(y^*)$. Hence, there is $u^*\in L^2(\Sigma)$ such that $y^{u^*}(T)\in[-\mu,\rho+\mu]$, a.e. in $\ooo^*$. This completes the proof of Theorem \ref{t3.1}.\hfill$\Box$

\begin{remark}\label{r3.4}\rm The controller $u^*$ as well as the corresponding mushy set $\ooo^{u^*,z}_T$ are dependent of $\ooo^*$ and so a realistic controlling process is that where the set $\ooo^*$ is sufficiently thin. Though the controllability theorem established above is not constructive, it leads to the theoretical conclusion that the mo\-ving interface of the two-Stefan problem is controllable. On the other hand, the optimal control problem \eqref{e3.8} with the state system \eqref{e2.1} provides an approximating solution for the controllability problem.
	\end{remark}

One could take the minimization problem \eqref{e3.8} with $z=y$ and the state system \eqref{e2.1} as an approximating control problem for \eqref{e3.1} though there are some technical problems regarding the convergence in this case of the solution~$u_\vp$.


\begin{thebibliography}{nn}
	
\bibitem{1} Atthey, D.R., A finite difference scheme for melting problems, {\it J. Inst. Math. Appl.}, 13 (1974), 353-366.

\bibitem{2a} Barbu, V., {\it Controllability and Stabilization of Parabolic Equations}, Birkh\"auser \& Springer International Publishing 2018.

\bibitem{2} Barbu, V., {\it Optimal Control of Variational Inequalities}, Pitman, London. Melbourne, 1983.

\bibitem{3} Barbu, V., {\it Analysis and Control of Nonliner Infinite Dimensional systems}, Academic Press, Boston. San Diego. New York, 1993.

\bibitem{3a} Barbu, V., Da Prato, G., Zolesio, J.P., Feedback controllability of the free boundary of the one phase Stefan problem, {\it Differential and Integral Equations}, 4 (1991), 225-239.

\bibitem{4} Barbu, V., Marinoschi, G., Existence for a time dependent rainfall infiltration model with a blowing up diffusivity, {\it Nonlinear Analysis. Real world Applications}, 5 (2004), 231-245.

\bibitem{5} Caginalp, G., Analysis of a phase field model of a free boundry, {\it Arch. Rat. Mech. Anal.}, 92 (1996), 206-245.

\bibitem{6} Koga, S., Krstic, M., Single boundary control of the two phase Stefan problem, {\it Systems \& Control Letters}, 135 (2020).

\bibitem{7} Fursikov, A., Imanuvilov, O., {\it Controllability of Evolution equations}, Lecture Notes, vol. 34, Seoul National University, 1996.

\bibitem{8} Hoffman, J., Sprekels, J., Real time control of free boundary in  a two-phase Stefan problem, {\it  Numer. Funct. Anal.}, 5 (1982), 47-76.

\bibitem{9} Lacely, H., Shillar, A., Mushy region in Stefan problem, {\it J. Inst. Math. Appl.}, 30 (1983), 303-313.

\bibitem{10} Lady\v zenskaja, O.A., Solonnikov, N.N., Uralceva, N.N., {\it Linear and Quasi-linear equations of Parabolic Type}, Translations of Mathematical Monographs, vol. 13, American Mathematical Society, 1988.

\bibitem{11} Lions, J.L., {\it Quelques m\'ethodes de resolution des probl\`mes aux limites non lin\'eaires}, Dunod-Gauthier-Villars, Paris, 1969.

\bibitem{12} Miranville, A., Morosanu, C., {\it Qualitative Analysis for the Ma\-the\-ma\-tical Models of Phase Separation and Transition. Applications}, Dif\-fe\-ren\-tial Equations \& Dynamical Systems, vol. 7, American Institute of Mathematical Sciences, 2020.

\bibitem{13} Marinoschi, G., {\it Functional Approach to Nonlinear Models of Water Flows in Solids}, Springer, 2006.

\bibitem{14} Saguez, C., Contr\^ole optimale de syst\`emes \`a fronti\`ere libre, Th\`ese l'Universit\'e de Technologie de Compi\`egne, 1980.

	
\end{thebibliography}
\end{document}